\documentclass[12pt]{article}

\newcommand{\pic}[2]{\BoxedEPSF{#1 scaled #2}}

\usepackage[hideboxes]{boxedeps}
\usepackage{a4,latexsym}
\newtheorem{corollary}{Corollary}

\newcommand\phm{\phantom{-}}

\newcommand\Bandij{\pic{Bandij.ART}{500}}
\newcommand\Cuspgraphic{\pic{Cuspgraphic.ART}{500}}
\newcommand\Doublesaddle{\pic{Doublesaddle.ART}{500}}
\newcommand\Discsaddle{\pic{Discsaddle.ART}{500}}
\newcommand\Labelledgraphic{\pic{Labelledgraphic.ART}{500}}
\newcommand\Saddlearcs{\pic{Saddlearcs.ART}{500}}
\newcommand\SixbraidA{\pic{SixbraidA.ART}{500}}
\newcommand\SixbraidB{\pic{SixbraidB.ART}{500}}
\newcommand\Stallingsfour{\pic{Stallingsfour.ART}{400}}
\newcommand\ThreebraidAB{\pic{ThreebraidAB.ART}{500}}
\newcommand\ThreebraidBA{\pic{ThreebraidBA.ART}{500}}
\newcommand\TwobraidAB{\pic{TwobraidAB.ART}{500}}
\newcommand\TwobraidBA{\pic{TwobraidBA.ART}{500}}

\newtheorem{theorem}{Theorem}

\newenvironment{Def}{\par\smallskip%
\noindent\textbf{Definition.}\  }%
{\par\smallskip}

\newenvironment{Proof}{\par\smallskip%
\noindent\textit{Proof : }\  }
{\hfill$\Box$\par\smallskip}
\begin{document}

\title{Mutual braiding and the band presentation of braid groups}

\author{H. R. Morton and M. Rampichini}
%
 


%

\maketitle

\begin{abstract}
This work is concerned with detecting when a closed braid and its axis are
`mutually braided' in the sense of Rudolph \cite{Rudolph2}. It deals with closed
braids which are fibred links, the simplest case being closed braids which
present the unknot. The geometric condition for mutual braiding refers to the
existence of  a close control on the way in which the whole family of fibre
surfaces meet the family of discs spanning the braid axis. We show how such a
braid can be presented naturally as a word in the `band generators' of the
braid group discussed by Birman, Ko and Lee \cite{BKL} in their recent account
of the band presentation of the braid groups. In this context we are able to
convert the conditions for mutual braiding into the existence of a   suitable
sequence of band relations and other moves on the braid word, and thus derive a
combinatorial method for deciding whether a braid is mutually braided.
\end{abstract} 

\section{Introduction} This is an account of  part of the second author's
doctoral dissertation \cite{Rampichini}. It is an
extension of the work of the first author on exchangeable braids dating from
about 1981, reported at the Sussex Low-dimensional Topology meeting of 1982,
\cite{Morton}. Among its antecedents are a problem of Stallings on representing
the unknot as a closed braid, and some constructions of fibred knots by
Goldsmith using closed braids. Subsequent work has been done by Rudolph
\cite{Rudolph2}, and by John Salkeld \cite{Salkeld} who showed that a braid
$\beta$ is exchangeable if and only if $\beta^2$ closes to a fibred link. A
more extended account has been written  \cite{Rampichini2}, giving further details of the
techniques and their applications to exchangeable braids.

\section{Exchangeable braids}
The original idea arose from the principle that many features of a braid
$\beta$ are best seen by looking at its closure $\hat\beta$ {\em along with its
axis} $A$ as a link $A\cup B$ with a distinguished unknotted component $A$, as
in figure 1.

\begin{center}\mbox{
\parbox[c]{.5\linewidth}{\begin{center}
\TwobraidAB\\

$\beta=\sigma_1$\end{center}}
\parbox[c]{.5\linewidth}{\begin{center}
\ThreebraidAB\\

$\beta=\sigma_1\sigma_2^{-1}$\end{center}}}
\\
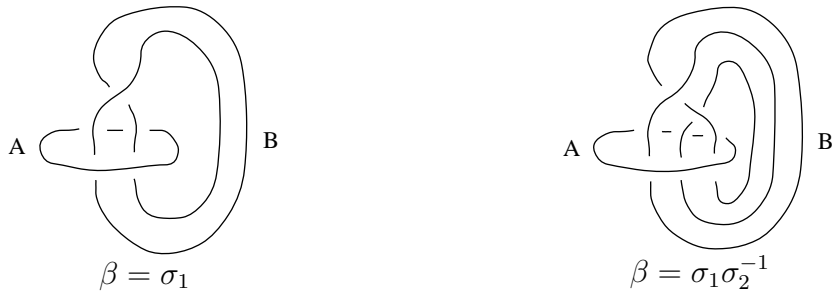
{Figure 1.}
\end{center}

Suppose conversely that we are given a link $A\cup B$ with $A$ unknotted. Then
there is a fibre projection $p_A:S^3-A\to S^1$ whose fibres
$A_\theta=p_A^{-1}(e^{i\theta})$ are discs.

\begin{Def} 
We say that $B$ is {\em braided rel} $A$ if it can be isotoped,
avoiding $A$, to meet all fibres transversely.
\end{Def}

Then $A\cup B$ can be viewed as a closed braid $B=\hat\beta$ with axis $A$; the
braid $\beta$ is determined up to conjugacy in $B_n$ where $n=|{\rm lk}(A,B)|$.
The link $A\cup B$ is thus a good way to capture the conjugacy class of a braid
$\beta$.

Exchangeable braids arise naturally when we look at such a link and ask whether
it is also possible for $B$ to be  the axis and $A$ the closed braid, in
other words whether we can have $B$ unknotted and $A=\hat\alpha$ braided rel
$B$.

\begin{Def}  We say that $A\cup B$ is {\em
exchangeably braided} if
$A$ and $B$ are unknotted, $B$ is braided rel $A$ and $A$ is braided rel $B$.
\end{Def}

The examples of figure 1 can be redrawn to see that both are exchangeably
braided, as in figure 2. In each of these cases there is even an isotopy of the
2-component link which interchanges $A$ and $B$, so that the braids $\alpha$ and
$\beta$ are the same up to conjugacy. 

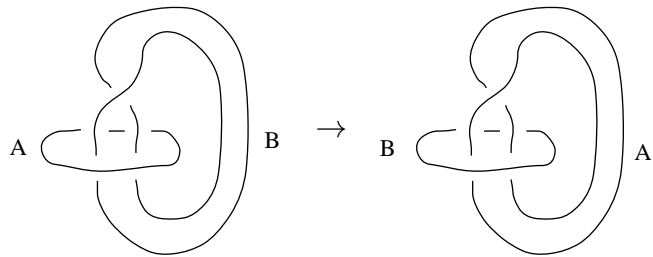
\begin{center}
\TwobraidAB\quad$\to$\quad\TwobraidBA \\
\ThreebraidAB\quad$\to$\quad\ThreebraidBA \\

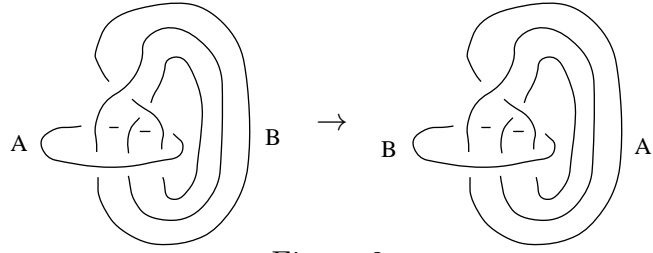
{Figure 2.}
\end{center}

This is not true in  general. For example the braid
$\beta$ shown in figure 3, which is built as a satellite of
$\sigma_1\sigma_2^{-1}$ with pattern $\sigma_1$, is exchangeable and the
exchanged braid $\alpha$, shown also in figure 3, is built as a satellite of
$\sigma_1$ with pattern $\sigma_1\sigma_2^{-1}$. However there can not be an
isotopy of the link $A\cup B$ interchanging the components, since its
2-variable Alexander polynomial $\Delta_{A\cup B}(a,b)=\sum c_{rs}a^rb^s$,
defined up to a power of $a$ and $b$, has coefficient  matrix
$$C=\pmatrix{0&0&1&0&0&0\cr
1&0&1&0&1&0\cr
0&0&1&0&0&0\cr
0&0&0&1&0&0\cr
0&1&0&1&0&1\cr
0&0&0&1&0&0\cr
}.$$ Now $C$ is not symmetric, and consequently 
$\alpha$ is not conjugate to
$\beta$.

\begin{center}
$\beta\quad=\quad$\SixbraidB \qquad$\alpha\quad=\quad$\SixbraidA\\
\medskip
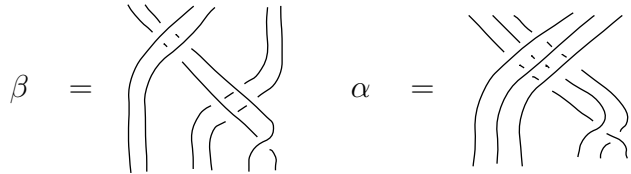
{Figure 3.}

\end{center}

\section{Bands}
The basic question of when a given braid $\beta$ is
exchangeable is most readily answered  in terms of 
the elementary $n$-braids $a_{ij}$ known as {\em embedded bands} in
\cite{Rudolph}, which we shall simply call {\em bands} in this paper.

\begin{Def}  The {\em band} $a_{ij}=a_{ji}$ is
the braid illustrated in figure 4, in which
the strings $i$ and $j$ interchange in the
positive sense in front of the other strings. 
\end{Def}

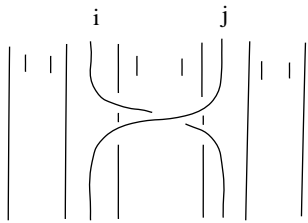
\begin{center}
\Bandij\\
\medskip
Figure 4.
\end{center}

\begin{Def} A {\em Stallings} $n$-braid is a
product of
$n-1$ bands or their inverses whose closure is connected. 
\end{Def}

The closure of such a braid is readily seen to be spanned by a disc
meeting the axis in $n$ points made up of one disc for each string, joined by
$n-1$ half-twisted bands. It is clearly
necessary for an exchangeable $n$-braid $\beta$ to have unknotted closure
$B=\hat\beta$, and for
$B$ to be spanned by a disc which meets the axis $A$ in $n$ points. This
property was shown in
\cite{Morton} to characterise Stallings braids up to conjugacy. It follows
 that any exchangeable braid is conjugate to a Stallings
braid.

 This condition alone is not enough, however; there has to be not just
one spanning disc, but a whole family $\{B_\varphi\}, \varphi\in[0,2\pi]$ which
all meet $A$ in $n$ points. Examples of non-exchangeable Stallings braids
are readily available for $n\ge4$. While the Stallings $4$-braids
$a_{13}(a_{23})^{-1}a_{24}$, illustrated in figure 5, and
$a_{13}a_{23}a_{24}$  are both exchangeable,  the Stallings braid
$\beta=a_{2 4}a_{23}a_{13}$ is {\em
not}. This can be detected  from the nature of the Alexander polynomial
$\Delta_{A\cup\hat\beta}(a,b)=\sum d_{rs}a^rb^s$ of any closed braid and axis. 
Such a polynomial always has the form $\det(aI-B(b))$ where $B(b)$ is the
$(n-1)\times(n-1)$ reduced Burau matrix of $\beta$, \cite{Morton}. The extreme
coefficients of powers of $a$ are then monomials in $b$ with coefficient
$\pm1$. A similar constraint must apply to the extreme coefficients of powers
of $b$ when the braid is exchangeable. In the case of
$\beta$ above, the matrix $D$ of coefficients of the Alexander polynomial is
$$D=\pmatrix{\phm1&\phm0&\phm0&\phm0\cr
\phm1&-1&\phm2&-1\cr-1&\phm2&-1&\phm1\cr\phm0&\phm0&\phm0&\phm1\cr};$$ although
the first and last rows have each a single entry
$\pm1$, this does not hold for the first and last columns, so in this case $A$
can not be braided rel
$\hat\beta$.

\begin{center}
\Stallingsfour\\
\medskip
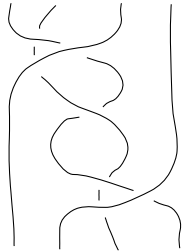
Figure 5.
\end{center}

 Birman, Ko and Lee \cite{BKL} have recently looked at the presentation of the
braid groups using bands as generators. The relations in this presentation
turn out to give a very good combinatorial method to help in determining
when a braid is exchangeable.

The relations themselves can be stated very simply as follows.
\begin{itemize}
\item  For $i>j>k$ we have
$a_{ij}a_{jk}=a_{jk}a_{ki}=a_{ki}a_{ij}$. 

The common product of
these pairs of bands can be visualised as a positive twist on the strings
$i,j$ and $k$ through one-third of a full turn. 

\item The bands $a_{p q}$ and $a_{rs}$ commute if $p,q,r$ and $s$ are
all different, and the pairs $(p,q)$ and $(r,s)$ do not interlace.

 Thus
$a_{14}a_{23}=a_{23}a_{14}$ but $a_{13}a_{24}\not
=a_{24}a_{13}$.
\end{itemize}

\section{Generalised exchangeable braids}

In trying to determine which Stallings braids are exchangeable
we follow Goldsmith \cite{Goldsmith} in using the term {\em generalised axis}
for an oriented fibred link $A$, not necessarily the unknot, and extending the
definition of a closed braid to include a curve $B$ which meets all the fibres
$A_\theta$ of the fibration $p_A:S^3-A\to S^1$ transversely. In this case we
also extend our previous terminology to say that $B$ is braided rel $A$ when
$A$ is fibred and the fibration can be chosen so that $B$ meets all the fibres
transversely, and thus meets each fibre in $n=|{\rm lk}(A,B)|$ points.

\begin{Def}  Call  $A\cup B$  a {\em
(generalised) exchangeable link} if
$A$ and
$B$ are both fibred, $B$ is braided rel $A$ {\em and} $A$ is braided rel
$B$.
\end{Def}

In the most general case, when both $A$ and $B$ are knotted, our analysis is
still rather sketchy, even in the case of braid index 1.

In what follows we shall restrict attention to the case in which one component,
$A$ say, is unknotted, when $B$ will determine a classical braid $\beta$ up to
conjugacy. We  call $\beta$ a {\em generalised exchangeable braid} if the
axis $A$ is braided rel $B=\hat\beta$ as above.

\begin{theorem} Every generalised exchangeable braid $\beta$ is
conjugate to a product of $k$ bands, such that
the resulting banded surface made up from one disc for each braid string joined
by $k$ half-twisted bands forms  a fibre surface for $B=\hat\beta$.
\end{theorem}

\begin{Proof}
We may suppose that  $\beta$ is an $n$-braid such that $B=\hat\beta$ is fibred
and whose axis $A$ is braided rel $B$. Then $A$ meets each fibre $B_\varphi$ in
$n$ points, with all intersections in the same sense.
Concentrate  on just one of these fibre surfaces $B_0$ and its intersection with
all the discs $A_\theta$ which span the axis $A$, or equivalently look at the  
function $p_A|B_0-A\to S^1$. We can assume that the intersections are transverse
except for finitely many values of $\theta$ where there is either a saddle or a
centre. Now $B_0$ has minimal genus among spanning surfaces of $B$ and all its
intersections with $A$ are in the same sense. Hence, following standard
arguments as in
\cite{Rudolph}, we can isotop it to eliminate centres. There are then
$n-\chi(B_0)$ saddles, each of which lies on a component of $A_\theta\cap B_0$
as in figure 6. 
\begin{center}
\Saddlearcs\\
\medskip
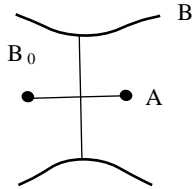
{Figure 6.}
\end{center}

The remainder of $B_0$ is foliated by the arcs of $A_\theta\cap B_0$ in a way
which is determined up to isotopy by the position of the saddles. Each saddle
 determines two crossing arcs in $B_0$, one which joins two  points of the
boundary $B$ and one which connects two of the $n$ points of $A$ lying in
$B_0$. Within $B_0$ take a neighbourhood $N$ of the points of $A$ and the cross
arcs of the $k=n-\chi(B_0)$ saddles which join them, whose boundary meets  the
foliation transversely. The
resulting subsurface differs from
$B_0$ only by a collar of the boundary $B$, so that $B$ can be isotoped
through
$B_0$, avoiding $A$, to form the boundary of $N$. We may choose a sufficiently
close neighbourhood $N$ as a surface in
$S^3$ made up of
$n$ disc neighbourhoods of the intersections with the axis together with
neighbourhoods of the $k$ cross-arcs each at a different level of $\theta$. The
close  neighbourhood of each such arc within $N$ lies as a half-twisted band
when viewed with $A$ as axis, while its boundary  still lies as a closed
braid. As a result,
$B_0$ is isotopic to a  banded surface $N$ whose boundary can replace the
original curve $B$ without loss. When the points of intersection of $B_0$ with
the axis are labelled $1,\ldots,n$ in order around $A$ then the neighbourhood
of a saddle arc joining points $i$ and $j$ becomes a band $(a_{ij})^{\pm 1}$,
whose sign depends on the direction of tangency of the surfaces $B_0$ and
$A_\theta$ at that saddle. The boundary can  be presented as the closure of a 
braid consisting of the product of these $k=n-\chi(B_0)$ bands in the order of
the value of
$\theta$ at each saddle. 
\end{Proof}

We may thus assume that any exchangeable braid is, up to conjugacy, a product
of bands and that a  fibre surface for its
closure is the resulting banded surface. In the case when $B$ is unknotted this
is the result from \cite{Morton} about Stallings braids. 

\begin{corollary} There are 
only a finite number of possible exchangeable links $A\cup B$ with unknotted
$A$, and a given genus for $B$ and linking number $n$. 
\end{corollary}

\begin{Proof}
The component $B$ can be written as the closure of an $n$-braid $\beta$ with
axis $A$. By the theorem we may assume that $\beta$ is the product of $k$
bands, where $k$ is determined by the genus of
$B$ and the linking number
$n$. There are only finitely many products of $k$ bands.
\end{Proof}

To  help decide exactly which products of bands {\em are}
exchangeable, we must look at the singular foliations induced on {\em every one}
of the family of fibre surfaces $B_\varphi$. At the same time we will see
singular foliations of each disc $A_\theta$ by the curves of $A_\theta\cap
B_\varphi$. 

We can best keep track of the foliations by plotting the values of
$(\theta,\varphi)\in S^1\times S^1$ for which $A_\theta$ and $B_\varphi$ do not
meet transversely. After a small isotopy we may assume that the singularities
of any intersections are  generic, so that the graphic of singular values
may look something like figure 7.

\begin{center}
\Cuspgraphic\\

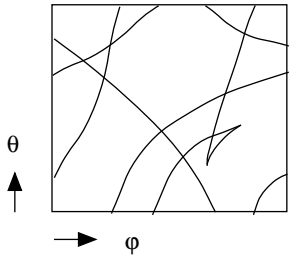
Figure 7.
\end{center}

At a general singular
point
$(\theta,\varphi)$ there is either a single saddle or centre tangency between
$A_\theta$ and
$B_\varphi$. For any $\varphi$ we can calculate the Euler characteristic of
 the $n$-punctured surface $B_{\varphi} - A$ as the number of centres minus the
number of saddles in the foliation of $B_\varphi$, which we can read from the
graphic on the line of the chosen value of
$\varphi$. Since this number must be constant ($=\chi(B_0)-n$) as
$\varphi$ varies we see that this sum remains fixed as the vertical line sweeps
through the graphic. Similarly each horizontal line meets the graphic so as to
contribute $\chi(A_\theta)-n=1-n$ to the Euler characteristic for the
$n$-punctured disc $A_\theta - B$.

Thus in the generic case no smooth local maxima or minima can occur either
vertically or horizontally on the lines of the graphic, which must consist of a
number of monotone increasing or decreasing lines with some simple crossings and
cusps. Each
monotone line represents either a persistent saddle or centre, by constancy of
the Euler characteristic. At a crossing, two tangencies generically occur at
distinct places and the saddles or centres continue through the crossings, while
at each cusp a line of centres and a line of saddles meet. The slope of a line
is positive when the surfaces
$A_\theta$ and
$B_\varphi$ have the same orientation at the point of tangency, and negative
when they have reverse orientations.

\section{Mutual braiding and labelled graphics}
To simplify further analysis we now impose Rudolph's condition that the two
fibrations form {\em mutually braided open books}, \cite{Rudolph2}.
Equivalently we require all the local tangencies between fibres of the two
families to be saddles, so that the graphic of singular values has no centres,
and hence no cusps. Rampichini proves in
\cite{Rampichini, Rampichini2} that, for an exchangeable link with one
unknotted component, this condition  can be assumed  without loss of generality.

\begin{Def}  A link $A\cup B$ with unknotted
$A$ which satisfies Rudolph's condition is
called {\em mutually braided}.
\end{Def}

The graphic of a mutually braided link then consists only of a number of
monotone lines, with simple crossings.  We can add some simple combinatorial
information to the graphic of a mutually braided link to describe the positions
of the saddles in each disc $A_\theta$. The boundary  of each disc
$A_\theta$ is the fixed curve $A$ which meets every fibre $B_\varphi$ in $n$
points. The $n$ points of $A$ where $\varphi=0$ as $1,\ldots,n$ dissect $A$ into
$n$ half-open segments. Label these
$1,\ldots,n$ in order around $A$.

 Where $A_\theta\cap B_\varphi$ contains a
single saddle tangency the cross-arc of the saddle  joins two points of $A$
lying in different segments $i$ and $j$. Label the corresponding point
$(\theta,\varphi)$ of the graphic of singular values with the pair $ij$. Along
any line of the graphic the label  remains constant   until we
reach a point where two lines of saddles cross, or we get to the transition
where
$\varphi$ increases through the level $\varphi=2\pi$ to return to $\varphi=0$.
At this transition the corresponding arcs of the saddles will move from segments
$i$ and $j$ to segments $i+1$ and $j+1$. 

Before describing the labelling of the graphic at a crossing it is worth
visualising the  singular
foliation of the disc $A_\theta$ for fixed $\theta$. A typical configuration
with $n=4$ is shown in figure 8. 

\begin{center}
\Discsaddle\\
\medskip
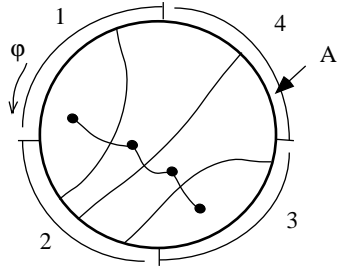
Figure 8.
\end{center}

The general line
$\theta={\rm constant}$ meets the graphic in $n-1$ points, each labelled by
some pair $ij$. These correspond to the cross arms of $n-1$ saddles, and the
labelling tells us in which segments of the circle $A$ their ends lie. They
divide the disc
$A_\theta$ into $n$ regions, each of which  contains one point of $B$; a
choice of this point and  the rest of the foliation is determined up to isotopy
by the saddle arms.  As $\theta$ increases the position of the saddles, and the
foliation changes, initially by an isotopy in which the end points of
the positive saddle arms rotate in one direction while those of the negative
arms rotate in the other. 

This evolution is viewed by Rampichini as a sort of cinematic film, in
which the screen is a disc with fixed boundary on which is projected the
foliation of $A_\theta$ with the parameter $\theta$ in the role of a  time
variable. The main part of the action is the movement of the saddle arms. This
is punctuated by isolated instants where two lines of saddles cross in the
graphic. In the film at each such critical time $\theta$ there are two saddles
having the same value of $\varphi$ and thus two singular points in the set
$A_\theta\cap B_\varphi$. These singular points either belong to two disjoint
components of $A_\theta\cap B_\varphi$, each consisting of simple saddle arms
connecting non-interlocking pairs of segments $p,q$ and $r,s$ say, or the two
saddles belong to the same component of $A_\theta\cap
B_\varphi$ which lies in $A_\theta$ as shown in figure 9. 

\begin{center}
\Doublesaddle\\

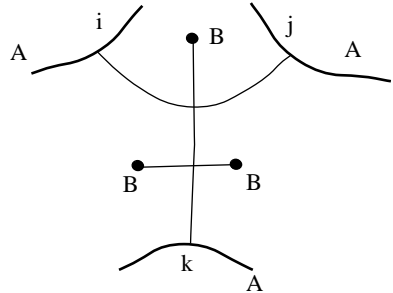
Figure 9.
\end{center}

Nothing very dramatic
happens in the film when disjoint saddle arms are involved; they simply
interchange levels of the  value $\varphi$, while in the graphic the labels
$pq$ and $rs$ continue through the crossing. Where the two saddles belong to
the same component, however, we see the saddle arms in the film coming together 
at one end and then reforming with one of the two arms remaining essentially the
same and the other changing radically to emerge from the other end of the
unchanged saddle arm. The whole process is determined by the `Christmas tree'
configuration of the intermediate position as in figure 9. The essential
information is the labelling of the three   segments
$i,j$ and $k$ of the end points, and the knowledge of which 
segment lies at the base of the tree. The labelling of the saddles on the
graphic before and  after the crossing depends on this combinatorial
information. As a summary of the film itself it is useful to think of simply
drawing the sequence of shots at the critical levels, with each Christmas tree
truncated to form a `T' and the signs of each saddle indicated so that the
direction of its evolution in the film is known. Rampichini presents such data,
with some extension, in what she terms a `film script', with these shots to
dictate the changes of scene. 

The labelling on the graphic is then quite constrained in its nature at the
crossing points. It is easier to describe the exact transitions of labelling at
the crossings when we change our viewpoint and consider what happens as we
change $\varphi$.

For a fixed $\varphi$ the graphic provides a complete list of all the saddles
in $B_\varphi$ which arise in the foliation by the intersections with the level
discs $A_\theta$, and the labelling shows which of the segments of $A$ are
connected by their cross-arms. The labelling on the graphic, read in order of
increasing $\theta$, determines a  word, $w_\varphi$, as a product of bands,
using 
$a_{ij}$ or its inverse according to the slope of the line on the graphic with
label
$ij$.

\begin{theorem} The band word $w_\varphi$ determined by a labelled graphic
arising from a mutually braided link    changes only
\begin{itemize}
\item when two lines of the graphic cross, or
\item when a line on the graphic passes between $\theta=0$ and
$\theta=2\pi$.
\end{itemize}
In the first case $w_\varphi$ changes by the application of a band relation on
passing the crossing. In the second case, which occurs exactly $n-1$
times, either a final band
$a_{ij}$  is transferred to the  beginning of the band word $w_\varphi$, or an
initial band
$a_{ij}^{-1}$ to the end, as $\varphi$ increases.
\end{theorem}

\begin{Proof} The first case is shown
by examining  the possible generic behaviour
close to the configuration in figure 9. Full
details  are given in \cite{Rampichini, Rampichini2}. The
second case is immediate from the monotonicity
of the lines and the fact that there are $n-1$
saddles at
$\theta=0$.
\end{Proof}

As a result, the labelled graphic arising from a mutally braided closed braid
$B$  determines a sequence of band words $w_\varphi$,  each representing $B$
and a spanning surface
$B_\varphi$. The sequence starts from the band word $w_0$, which
can be regarded as a given band word for the proposed exchangeable braid $B$
along with its spanning surface, and finishes with the word $w_{2\pi}$, which is
the original word $w_0$ with all its indices $i,j$ reduced by $1$ mod
$n$. They are connected by  a sequence of moves consisting of  an unspecified number of
applications of the band relations along with
 exactly
$n-1$ cyclings as specified in the theorem. 

By way of example, a labelled graphic for the exchangeable 4-braid with band
word $a_{34}(a_{12})^{-1}a_{23}$ is shown in figure 10.

\begin{center}
\Labelledgraphic\\

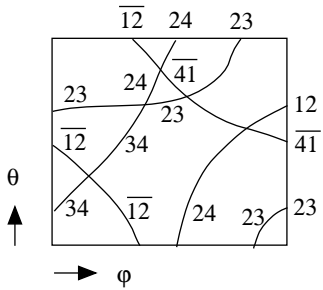
Figure 10.
\end{center}

The resulting sequence of band words, given simply by their indices, is listed
below.
\begin{center}
$(34)\overline{(12)}(23)\to \overline{(12)}(34)(23)\to
(34)(23)\overline{(12)}$\\
$\to (23)(24)\overline{(12)}\to
(23)\overline{(14)}(24)
\to(24)(23)\overline{(14)}$\\
$\to(24)\overline{(14)}(23)\to(23)(24)\overline{(14)}\to(23)\overline{(41)}(12).$
\end{center}

\section{Band algorithms for mutual braiding}
The previous section shows that if a braid determines a mutually braided open
book then there must be a labelled graphic giving rise to a suitably restricted
sequence of band words. It can be shown, with a little more argument, that 
there are only finitely many possible sequences which need to be considered when
asking if a given band word determines a mutually braided link, and if none of
them fulfil all the necessary conditions for arising from a labelled graphic
then the link is definitely not mutually braided.

The key to the converse of this result, which also covers Rampichini's extension
to exchangeable braids, lies in the reconstruction of a family of spanning
surfaces from an `admissibly' labelled graphic, \cite{Rampichini,Rampichini2}. 

\begin{Def}  A graphic of
monotone increasing and decreasing lines  in $S^1\times S^1$,
labelled by pairs of numbers $i,j$ between $1$ and $n$, is {\em admissibly
labelled} if
\begin{itemize}
\item   the labelled lines cross one basic circle $\theta=0$ in
$n-1$ points,
\item  the band words $w_\varphi$   read from the graphic  following the 
basic circles
$\varphi={\rm constant}$ of the other family satisfy the band relations at
every crossing,
\item 
the indices in the bands all change by 1
at the reference circle
$\varphi =0$,
\item  the $n-1$ labels on the
circle $\theta=0$ contain no interlocking pairs.
\end{itemize}
\end{Def}

\begin{theorem} Every admissibly labelled graphic arises from some exchangeable
braid.
\end{theorem}

\begin{Proof} Use the last condition in the
definition to draw a set of saddles in the
disc $A_0$ which correspond with the labelling
when $\theta=0$. The evolution of the graphic
as
$\theta$ increases then provides a film script, which is converted into a film
by filling in suitably the non-singular parts of the foliation at each level.
Use of standard models around the singular levels where there are crossings in
the graphic, coupled with  the monotonicity of the movement of the
saddle arms with increasing $\theta$, ensures that the family of 
surfaces in the solid torus generated in this way are non-singular everywhere,
\cite{Rampichini, Rampichini2}.
\end{Proof}

Consequently an admissibly labelled graphic with a given band word at
$\varphi=0$ will give rise to a mutually braided fibration for the
corresponding closed braid.

Rampichini also shows how the graphic for a general exchangeable braid, which
includes cusps and lines of centres, can be used to construct an admissibly
labelled graphic, and hence deduce that any exchangeable braid is mutually
braided
\cite{Rampichini,Rampichini2}.

 We thus
 have an algorithm to decide whether a band word $w_0$ determines a
mutually braided (and indeed an exchangeable) braid.
\medskip

\noindent{\bf Algorithm.}\quad
Look for a chain of moves to convert the word $w_0$ to the same word but with
all indices reduced by $1$ mod $n$. The moves must consist of a number of band
relations and exactly $n-1$ cyclings in which a positive band is moved from back
to front or a negative band from front to back of the word. It can be shown
that no repetitions of words between cyclings need occur in an admissibly
labelled graphic, so there are only a finite number of possible chains to
consider \cite{Rampichini,Rampichini2}.

If no chains satisfy these conditions then the braid is {\em not} mutually
braided. 

Given a chain satisfying these conditions try to realise it by an
admissibly labelled graphic. It may not be possible to do this, bearing in mind
that the lines in the graphic, which are determined combinatorially by the
chain as
$\varphi$ increases, are required to move strictly monotonically in $\theta$,
and must end at the same level that they begin. If any chain can be realised
then the braid {\em is} mutually braided, otherwise again it is not.
\medskip

As an example we apply this method to the Stallings braid
$w=a_{24}a_{23}a_{13}$ mentioned earlier to show that its closure is not
exchangeable. 

We must consider possible chains of moves which lead to the word
$w'=a_{13}a_{12}a_{24}$, and include exactly $3$ cyclings. Now no band
relations are possible in
$w$ so the first move in any chain must be a cycling, to get
$a_{13}a_{24}a_{23}$. Again no relations are possible, so we must cycle to
$a_{23}a_{13}a_{24}$. Once more there are no relations and we must cycle,
returning to $w$. Since we have used all $3$ cyclings there is nothing more we
can do; as we have not reached $w'$  we conclude that the braid is not
exchangeable.

In contrast we {\em can} find a chain and an admissibly labelled graphic
leading from the word
$v=a_{13}a_{23}a_{24}$ to $v'=a_{24}a_{12}a_{13}$, by making use of  band
relations such as $a_{23}a_{24}=a_{24}a_{34}$.

\section*{Acknowledgments} The  second author
was partially supported by a British
Council Fellowship Award during the academic
year 1996-97.

{\noindent Department of Mathematical Sciences\\
 University of Liverpool\\
 Liverpool L69 3BX\\ England.\\
morton@liv.ac.uk\\
http://www.liv.ac.uk/\~{ }su14\\{ }\\}

{\noindent Dipartimento di Matematica\\Universit\` a degli Studi\\ via
Saldini 50\\20133 Milano, Italy.\\
rampichini@mat.unimi.it}

\end{document}